\RequirePackage{fix-cm}
\documentclass[smallextended]{svjour3}       
\smartqed  
\usepackage{amsmath}
\allowdisplaybreaks
 \spnewtheorem{thm}{Theorem}{\bfseries}{\itshape}
 \spnewtheorem{cor}[thm]{Corollary}{\bfseries}{\itshape}

\usepackage{graphicx,amssymb}
\global\long\def\relphantom#1{\mathrel{\phantom{{#1}}}}
\global\long\def\Zp{\mathbb{Z}_{p}}

\begin{document}

\title{Symmetric identities of higher-order degenerate Euler polynomials
}


\author{Dae San Kim         \and
        Taekyun Kim 
}


\institute{Dae San Kim \at
              Department of Mathematics, Sogang University, Seoul 121-742, Republic
of Korea \\
              Tel.: +82-2-7058871\\
              Fax: +82-2-7058871\\
\email{dskim@sogang.ac.kr}
           \and
           Taekyun Kim \at
              Department of Mathematics, Kwangwoon University, Seoul 139-701, Republic
of Korea\\
\email{tkkim@kw.ac.kr}
}

\date{Received: date / Accepted: date}

\maketitle

\begin{abstract}
The purpose of this paper is to give some symmetric identities of
higher-order degenerate Euler polynomials derived from the
symmetric properties of the multivariate $p$-adic fermionic integrals
on $\mathbb{Z}_{p}$.
\end{abstract}

\section{Introduction}

Let $p$ be a fixed prime such that $p\equiv1\pmod{2}$. Throughout
this paper, $\Zp$, $\mathbb{Q}_{p}$, and $\mathbb{C}_{p}$ will denote
the ring of $p$-adic integers, the field of $p$-adic numbers and
the completion of the algebraic closure of $\mathbb{Q}_{p}$, respectively.
Let $\nu_{p}$ be the normalized exponential valuation of $\mathbb{C}_{p}$
with $\left|p\right|_{p}=p^{-\nu_{p}\left(p\right)}=p^{-1}$. 

Let $f\left(x\right)$ be a continuous function on $\Zp$. Then $p$-adic
fermionic integral on $\Zp$ is defined by Kim as
\begin{align}
\int_{\Zp}f\left(x\right)d\mu_{-1}\left(x\right) & =\lim_{N\rightarrow\infty}\sum_{x=0}^{p^{N}-1}f\left(x\right)\mu_{-1}\left(x+p^{N}\mathbb{Z}_{p}\right)\label{eq:1}\\
 & =\lim_{N\rightarrow\infty}\sum_{x=0}^{p^{N}-1}f\left(x\right)\left(-1\right)^{x},\quad\left(\text{see \cite{key-8}}\right).\nonumber 
\end{align}

From (\ref{eq:1}), we note that 
\begin{equation}
\int_{\Zp}f\left(x+1\right)d\mu_{-1}\left(x\right)+\int_{\Zp}f\left(x\right)d\mu_{-1}\left(x\right)=2f\left(0\right),\label{eq:2}
\end{equation}
and 
\begin{align}
&\int_{\Zp}f\left(x+n\right)d\mu_{-1}\left(x\right)+\left(-1\right)^{n-1}\int_{\Zp}f\left(x\right)d\mu_{-1}\left(x\right)\nonumber\\
&=2\sum_{l=0}^{n-1}\left(-1\right)^{n-1-l}f\left(l\right),\quad\left(\text{see \cite{key-4}}\right).\label{eq:3}
\end{align}

As is well known, the Euler polynomials are defined by the generating
function 
\begin{equation}
\frac{2}{e^{t}+1}e^{xt}=\sum_{n=0}^{\infty}E_{n}\left(x\right)\frac{t^{n}}{n!},\quad\left(\text{see \cite{key-4,key-8}}\right).\label{eq:4}
\end{equation}

When $x=0$, $E_{n}=E_{n}\left(0\right)$ are called the Euler numbers. 

For $r\in\mathbb{N}$, the higher-order Euler polynomials are given
by 
\begin{equation}
\left(\frac{2}{e^{t}+1}\right)^{r}e^{xt}=\sum_{n=0}^{\infty}E_{n}^{\left(r\right)}\left(x\right)\frac{t^{n}}{n!},\quad\left(\text{see \cite{key-1,key-2,key-3,key-4,key-5,key-6,key-7,key-9,key-10}}\right)\label{eq:5}
\end{equation}

In particular, $x=0$, $E_{n}^{\left(r\right)}=E_{n}^{\left(r\right)}\left(0\right)$
are called the Euler numbers of order $r$. 

From (\ref{eq:2}), we can easily derive the following equation: 
\begin{align}
 & \int_{\Zp}\cdots\int_{\Zp}e^{\left(x_{1}+\cdots+x_{r}+x\right)t}d\mu_{-1}\left(x_{1}\right)\cdots d\mu_{-1}\left(x_{r}\right)\label{eq:6}\\
 & =\left(\frac{2}{e^{t}+1}\right)^{r}e^{xt}\nonumber \\
 & =\sum_{n=0}^{\infty}E_{n}^{\left(r\right)}\left(x\right)\frac{t^{n}}{n!}.\nonumber 
\end{align}

Thus, by (\ref{eq:6}), we get 
\begin{equation}
\int_{\Zp}\cdots\int_{\Zp}\left(x_{1}+\cdots+x_{r}+x\right)^{n}d\mu_{-1}\left(x_{1}\right)\cdots d\mu_{-1}\left(x_{r}\right)=E_{n}^{\left(r\right)}\left(x\right),\quad\left(n\ge0\right).\label{eq:7}
\end{equation}

Carlitz introduced the degenerate Euler polynomials given
by the generating function 
\begin{equation}
\left(\frac{2}{\left(1+\lambda t\right)^{\frac{1}{\lambda}}+1}\right)\left(1+\lambda t\right)^{\frac{x}{\lambda}}=\sum_{n=0}^{\infty}\mathcal{E}_{n}\left(x\mid\lambda\right)\frac{t^{n}}{n!},\quad\left(\text{see \cite{key-5}}\right).\label{eq:7-1}
\end{equation}
When $x=0$, $\mathcal{E}_{n}\left(0\mid\lambda\right)=\mathcal{E}_{n}\left(\lambda\right)$
are the degenerate Euler numbers. Note that $\lim_{\lambda\rightarrow0}\mathcal{E}_{n,\lambda}\left(x\mid\lambda\right)=E_{n}\left(x\right).$ 

For $r\in\mathbb{N}$, the higher-order degenerate Euler polynomials
are also given by the generating function  
\begin{equation}
\left(\frac{2}{\left(1+\lambda t\right)^{\frac{1}{\lambda}}+1}\right)^{r}\left(1+\lambda t\right)^{\frac{x}{\lambda}}=\sum_{n=0}^{\infty}\mathcal{E}_{n}^{\left(r\right)}\left(x\mid\lambda\right)\frac{t^{n}}{n!}.\label{eq:8}
\end{equation}

When $x=0$, $\mathcal{E}_{n}^{\left(r\right)}\left(0\mid\lambda\right)=\mathcal{E}_{n}^{\left(r\right)}\left(\lambda\right)$
are called the higher-order degenerate Euler numbers (see \cite{key-2}). 

In \cite{key-8}, Kim and Kim showed that the degenerate Euler polynomials
can be represented by a $p$-adic integral on $\Zp$. Recently,
several researchers have studied the symmetric identities of higher-order
Euler polynomials derived from the symmetric properties
of $p$-adic integrals on $\Zp$(see \cite{key-1,key-2,key-3,key-4,key-5,key-6,key-7,key-8,key-9,key-10,key-11,key-12,key-13,key-14,key-15,key-16,key-17,key-18,key-19}). 
 
In this paper, we investigate some properties of symmetry for the
multivariate $p$-adic fermionic integrals on $\Zp$. From our investigation,
we derive some identities of symmetry for the higher-order degenerate
Euler polynomials.

\section{Identities of symmetry for the higher-order degenerate Euler polynomials}

In this section, we assume that $\lambda,t\in\mathbb{C}_{p}$ with $\left|\lambda t\right|_{p}<p^{-\frac{1}{p-1}}$. From (\ref{eq:2}),
we can derive the following equation:
\begin{align}
 & \int_{\Zp}\cdots\int_{\Zp}\left(1+\lambda t\right)^{\frac{x_{1}+\cdots+x_{r}+x}{\lambda}}d\mu_{-1}\left(x_{1}\right)\cdots d\mu_{-1}\left(x_{r}\right)\label{eq:9}\\ 
 & =\left(\frac{2}{\left(1+\lambda t\right)^{\frac{1}{\lambda}}+1}\right)^{r}\left(1+\lambda t\right)^{\frac{x}{\lambda}}\nonumber \\
 & =\sum_{n=0}^{\infty}\mathcal{E}_{n}^{\left(r\right)}\left(x\mid\lambda\right)\frac{t^{n}}{n!}.\nonumber 
\end{align}

Thus, by (\ref{eq:9}), we get, for $n\geq0$,
\begin{equation}
\int_{\Zp}\cdots\int_{\Zp}\left(x_{1}+\cdots+x_{r}+x\mid\lambda\right)_{n}d\mu_{-1}\left(x_{1}\right)\cdots d\mu_{-1}\left(x_{r}\right)=\mathcal{E}_{n}^{\left(r\right)}\left(x\mid\lambda\right),\label{eq:10}
\end{equation}
where 
\begin{equation}
\left(x\mid\lambda\right)_{n}=x\left(x-\lambda\right)\cdots\left(x-\lambda\left(n-1\right)\right),\quad\left(n\ge0\right).\label{eq:11}
\end{equation}

From (\ref{eq:11}), we note that 
\begin{equation}
\left(x\mid\lambda\right)_{n}=\sum_{l=0}^{n}S_{1}\left(n,l\right)\lambda^{n-l}x^{l},\quad\left(n\ge0\right),\label{eq:12}
\end{equation}
where $S_{1}\left(n,l\right)$ is the Stirling number of the first
kind. 

From (\ref{eq:3}), we have, for $m\geq0$,
\begin{align}
&2\sum_{l=0}^{n-1}\left(-1\right)^{n-1-l}\left(l\mid\lambda\right)_{m}\label{eq:13}\\
& =\int_{\Zp}\left(x+n\mid\lambda\right)_{m}d\mu_{-1}\left(x\right)+\left(-1\right)^{n-1}\int_{\Zp}\left(x\mid\lambda\right)_{m}d\mu_{-1}\left(x\right).\nonumber
\end{align}

For $n\in\mathbb{N}$ with $n\equiv1\pmod{2}$, we get 
\begin{equation}
2\sum_{l=0}^{n-1}\left(-1\right)^{l}\left(l\mid\lambda\right)_{m}=\mathcal{E}_{m}\left(n\mid\lambda\right)+\mathcal{E}_{m}\left(\lambda\right),\quad\left(m\ge0\right).\label{eq:14}
\end{equation}

Let us define $\tilde{S}_{k}\left(n\mid\lambda\right)$ as follows:
\begin{equation}
\tilde{S}_{k}\left(n\mid\lambda\right)=\sum_{l=0}^{n}\left(-1\right)^{l}\left(l\mid\lambda\right)_{k},\quad\left(k,n\ge0\right).\label{eq:15}
\end{equation}

Then, we note that 
\begin{equation}
\lim_{\lambda\rightarrow0}\tilde{S_{k}}\left(n\mid\lambda\right)=\sum_{l=0}^{n}\left(-1\right)^{l}l^{k}=\tilde{S_{k}}\left(n\right).\label{eq:16}
\end{equation}

Let $n\in\mathbb{N}$ with $n\equiv1\pmod{2}$. Then we get 
\begin{align}
 & \frac{2\int_{\Zp}\left(1+\lambda t\right)^{\frac{x}{\lambda}}d\mu_{-1}\left(x\right)}{\int_{\Zp}\left(1+\lambda t\right)^{\frac{nx}{\lambda}}d\mu_{-1}\left(x\right)}\label{eq:17}\\
 & =\int_{\Zp}\left(1+\lambda t\right)^{\frac{n+x}{\lambda}}d\mu_{-1}\left(x\right)+\int_{\Zp}\left(1+\lambda t\right)^{\frac{x}{\lambda}}d\mu_{-1}\left(x\right)\nonumber \\
 & =2\sum_{l=0}^{n-1}\left(-1\right)^{l}\left(1+\lambda t\right)^{\frac{l}{\lambda}}\nonumber \\
 & =2\sum_{k=0}^{\infty}\tilde{S_{k}}\left(n-1\mid\lambda\right)\frac{t^{k}}{k!}.\nonumber 
\end{align}

Let $w_{1},w_{2}\in\mathbb{N}$, with $w_{1}\equiv1\pmod{2}$,
$w_{2}\equiv1\pmod{2}$.
For $m\in\mathbb{N}$, we define 
\begin{align}
 & K^{\left(m\right)}\left(w_{1},w_{2}\mid\lambda\right)\label{eq:18}\\
 & =\left(\frac{2}{\left(1+\lambda t\right)^{\frac{w_{1}}{\lambda}}+1}\right)^{m}\left(1+\lambda t\right)^{\frac{w_{1}w_{2}x}{\lambda}}\left(\left(1+\lambda t\right)^{\frac{w_{1}w_{2}}{\lambda}}+1\right)\nonumber \\
 & \relphantom =\times\left(\frac{2}{\left(1+\lambda t\right)^{\frac{w_{2}}{\lambda}}+1}\right)^{m}\frac{1}{2}\left(1+\lambda t\right)^{\frac{w_{1}w_{2}}{\lambda}y}.\nonumber 
\end{align}

From (\ref{eq:3}), we note that 
\begin{align}
K^{\left(m\right)}\left(w_{1},w_{2}\mid\lambda\right) & =\frac{\int_{\Zp^{m}}\left(1+\lambda t\right)^{\frac{w_{1}}{\lambda}\left(x_{1}+\cdots+x_{m}+w_{2}x\right)}d\mu_{-1}\left(x_{1}\right)\cdots d\mu_{-1}\left(x_{m}\right)}{\int_{\Zp}\left(1+\lambda t\right)^{\frac{w_{1}w_{2}}{\lambda}x}d\mu_{-1}\left(x\right)}\label{eq:19}\\
 & \relphantom =\times\int_{\Zp^{m}}\left(1+\lambda t\right)^{\frac{w_{2}}{\lambda}\left(x_{1}+\cdots+x_{m}+w_{1}y\right)}d\mu_{-1}\left(x_{1}\right)\cdots d\mu_{-1}\left(x_{m}\right),\nonumber 
\end{align}
where 
\begin{align}
&\int_{\Zp^{m}}f\left(x_{1},\dots,x_{m}\right)d\mu_{-1}\left(x_{1}\right)\cdots d\mu_{-1}\left(x_{m}\right)\label{eq:20}\\
&=\int_{\Zp}\cdots\int_{\Zp}f\left(x_{1},\dots,x_{m}\right)d\mu_{-1}\left(x_{1}\right)\cdots d\mu_{-1}\left(x_{m}\right).\nonumber
\end{align}

It is easy to see that $K^{\left(m\right)}\left(w_{1},w_{2}\mid\lambda\right)$
is symmetric in $w_{1}$ and $w_{2}$. Now, we observe that 
\begin{align}
 & K^{\left(m\right)}\left(w_{1},w_{2}\mid\lambda\right)\label{eq:21}\\
 & =\left(1+\lambda t\right)^{\frac{w_{1}w_{2}}{\lambda}x}\int_{\Zp^{m}}\left(1+\lambda t\right)^{\frac{w_{1}}{\lambda}\left(x_{1}+\cdots+x_{m}\right)}d\mu_{-1}\left(x_{1}\right)\cdots d\mu_{-1}\left(x_{m}\right)\nonumber \\
 & \relphantom =\times\frac{\int_{\Zp}\left(1+\lambda t\right)^{\frac{w_{2}}{\lambda}x_{m}}d\mu_{-1}\left(x_{m}\right)}{\int_{\Zp}\left(1+\lambda t\right)^{\frac{1}{\lambda}w_{1}w_{2}x}d\mu_{-1}\left(x\right)}\nonumber \\
 & \relphantom =\times\left(1+\lambda t\right)^{\frac{1}{\lambda}w_{1}w_{2}y}\int_{\Zp^{m-1}}\left(1+\lambda t\right)^{\frac{w_{2}}{\lambda}\left(x_{1}+\cdots+x_{m-1}\right)}d\mu_{-1}\left(x_{1}\right)\cdots d\mu_{-1}\left(x_{m-1}\right).\nonumber 
\end{align}

By (\ref{eq:9}), we get 
\begin{align}
 & \left(1+\lambda t\right)^{\frac{w_{1}w_{2}}{\lambda}x}\int_{\Zp^{m}}\left(1+\lambda t\right)^{\frac{w_{1}}{\lambda}\left(x_{1}+\cdots+x_{m}\right)}d\mu_{-1}\left(x_{1}\right)\cdots d\mu_{-1}\left(x_{m}\right)\label{eq:22}\\
 & =\left(\frac{2}{\left(1+\lambda t\right)^{\frac{w_{1}}{\lambda}}+1}\right)^{m}\left(1+\lambda t\right)^{\frac{w_{1}w_{2}}{\lambda}x}\nonumber \\
 & =\sum_{n=0}^{\infty}\mathcal{E}_{n}^{\left(m\right)}\left(w_{2}x\left|\frac{\lambda}{w_{1}}\right.\right)w_{1}^{n}\frac{t^{n}}{n!}.\nonumber 
\end{align}

From (\ref{eq:17}), (\ref{eq:21}) and (\ref{eq:22}), we have 
\begin{align}
 & K^{\left(m\right)}\left(w_{1},w_{2,}\mid\lambda\right)\label{eq:23}\\
 & =\sum_{l=0}^{\infty}\mathcal{E}_{l}^{\left(m\right)}\left(w_{2}x\left|\frac{\lambda}{w_{1}}\right.\right)w_{1}^{l}\frac{t^{l}}{l!}\sum_{k=0}^{\infty}\tilde{S_{k}}\left(w_{1}-1\left|\frac{\lambda}{w_{2}}\right.\right)\frac{w_{2}^{k}}{k!}t^{k}\nonumber \\
 & \relphantom =\times\sum_{i=0}^{\infty}\mathcal{E}_{i}^{\left(m-1\right)}\left(w_{1}y\left|\frac{\lambda}{w_{2}}\right.\right)\frac{w_{2}^{i}}{i!}t^{i}\nonumber \\
 & =\sum_{l=0}^{\infty}\mathcal{E}_{l}^{\left(m\right)}\left(w_{2}x\left|\frac{\lambda}{w_{1}}\right.\right)w_{1}^{l}\frac{t^{l}}{l!}\nonumber\\
 & \relphantom =\times\sum_{j=0}^{\infty}\sum_{k=0}^{j}\tilde{S_{k}}\left(w_{1}-1\left|\frac{\lambda}{w_{2}}\right.\right)\binom{j}{k}w_{2}^{k}w_{2}^{j-k}\mathcal{E}_{j-k}^{\left(m-1\right)}\left(w_{1}y\left|\frac{\lambda}{w_{2}}\right.\right)\frac{t^{j}}{j!}\nonumber \\
 & =\sum_{n=0}^{\infty}\sum_{j=0}^{n}\binom{n}{j}w_{2}^{j}w_{1}^{n-j}\mathcal{E}_{n-j}^{\left(m\right)}\left(w_{2}x\left|\frac{\lambda}{w_{1}}\right.\right)\nonumber\\
 & \relphantom =\times\sum_{k=0}^{j}\tilde{S_{k}}\left(w_{1}-1\left|\frac{\lambda}{w_{2}}\right.\right)\binom{j}{k}\mathcal{E}_{j-k}^{\left(m-1\right)}\left(w_{1}y\left|\frac{\lambda}{w_{2}}\right.\right)\frac{t^{n}}{n!}.\nonumber 
\end{align}

On the other hand, 
\begin{align}
 & K^{\left(m\right)}\left(w_{1},w_{2}\mid\lambda\right)\label{eq:24}\\
 & =\sum_{n=0}^{\infty}\sum_{j=0}^{n}\binom{n}{j}w_{1}^{j}w_{2}^{n-j}\mathcal{E}_{n-j}^{\left(m\right)}\left(w_{1}x\left|\frac{\lambda}{w_{2}}\right.\right)\nonumber\\
 & \relphantom =\times\sum_{k=0}^{j}\tilde{S_{k}}\left(w_{2}-1\left|\frac{\lambda}{w_{1}}\right.\right)\binom{j}{k}\mathcal{E}_{j-k}^{\left(m-1\right)}\left(w_{2}y\left|\frac{\lambda}{w_{1}}\right.\right)\frac{t^{n}}{n!}.\nonumber 
\end{align}

Therefore, by (\ref{eq:23}) and (\ref{eq:24}), we obtain the following
theorem.
\begin{thm}
\label{thm:1} For $w_{1},w_{2}\in\mathbb{N}$, with $w_{1}\equiv1\pmod{2}$,
$w_{2}\equiv1\pmod{2}$, $n\ge0$ and $m\in\mathbb{N}$, we have 
\begin{align*}
 & \sum_{j=0}^{n}\binom{n}{j}w_{2}^{j}w_{1}^{n-j}\mathcal{E}_{n-j}^{\left(m\right)}\left(w_{2}x\left|\frac{\lambda}{w_{1}}\right.\right)\sum_{k=0}^{j}\tilde{S_{k}}\left(w_{1}-1\left|\frac{\lambda}{w_{2}}\right.\right)\binom{j}{k}\mathcal{E}_{j-k}^{\left(m-1\right)}\left(w_{1}y\left|\frac{\lambda}{w_{2}}\right.\right)\\
 & =\sum_{j=0}^{n}\binom{n}{j}w_{1}^{j}w_{2}^{n-j}\mathcal{E}_{n-j}^{\left(m\right)}\left(w_{1}x\left|\frac{\lambda}{w_{2}}\right.\right)\sum_{k=0}^{j}\tilde{S_{k}}\left(w_{2}-1\left|\frac{\lambda}{w_{1}}\right.\right)\binom{j}{k}\mathcal{E}_{j-k}^{\left(m-1\right)}\left(w_{2}y\left|\frac{\lambda}{w_{1}}\right.\right).
\end{align*}

\end{thm}
Let $y=0$ and $m=1$ in Theorem \ref{thm:1}. Then we have the following
theorem.
\begin{thm}
\label{thm:2} For $n\ge0$, $w_{1},w_{2}\in\mathbb{N}$ with $w_{1}\equiv1,w_{2}\equiv1\pmod{2}$,
we have 
\begin{align*}
 & \sum_{j=0}^{n}\binom{n}{j}w_{2}^{j}w_{1}^{n-j}\mathcal{E}_{n-j}\left(w_{2}x\left|\frac{\lambda}{w_{1}}\right.\right)\tilde{S_{j}}\left(w_{1}-1\left|\frac{\lambda}{w_{2}}\right.\right)\\
 & =\sum_{j=0}^{n}\binom{n}{j}w_{1}^{j}w_{2}^{n-j}\mathcal{E}_{n-j}\left(w_{1}x\left|\frac{\lambda}{w_{2}}\right.\right)\tilde{S_{j}}\left(w_{2}-1\left|\frac{\lambda}{w_{1}}\right.\right).
\end{align*}

\end{thm}
In particular, if we take $w_{2}=1$ in Theorem \ref{thm:2}, then
we obtain the following corollary.
\begin{cor}
\label{cor:3} For $w_{1}\in\mathbb{N}$ with $w_{1}\equiv1\pmod{2}$,
$n\ge0$, we have 
\[
\mathcal{E}_{n}\left(w_{1}x\mid\lambda\right)=\sum_{j=0}^{n}\binom{n}{j}w_{1}^{n-j}\mathcal{E}_{n-j}\left(x\left|\frac{\lambda}{w_{1}}\right.\right)\tilde{S_{j}}\left(w_{1}-1\mid\lambda\right).
\]
 
\end{cor}
From (\ref{eq:18}), we have 
\begin{align}
 & K^{\left(m\right)}\left(w_{1},w_{2}\mid\lambda\right)\label{eq:25}\\
 & =\left(\frac{2}{\left(1+\lambda t\right)^{\frac{w_{1}}{\lambda}}+1}\right)^{m}\left(1+\lambda t\right)^{\frac{w_{1}w_{2}}{\lambda}x}\nonumber \\
 & \relphantom =\times\frac{\left(1+\lambda t\right)^{\frac{w_{1}w_{2}}{\lambda}}+1}{\left(1+\lambda t\right)^{\frac{w_{2}}{\lambda}}+1}\left(\frac{2}{\left(1+\lambda t\right)^{\frac{w_{2}}{\lambda}}+1}\right)^{m-1}\left(1+\lambda t\right)^{\frac{w_{1}w_{2}}{\lambda}y}\nonumber \\
 & =\left(\frac{2}{\left(1+\lambda t\right)^{\frac{w_{1}}{\lambda}}+1}\right)^{m}\sum_{i=0}^{w_{1}-1}\left(-1\right)^{i}\left(1+\lambda t\right)^{\frac{w_{1}w_{2}}{\lambda}x+\frac{w_{2}}{\lambda}i}\nonumber \\
 & \relphantom =\times\left(\frac{2}{\left(1+\lambda t\right)^{\frac{w_{2}}{\lambda}}+1}\right)^{m-1}\left(1+\lambda t\right)^{\frac{w_{1}w_{2}}{\lambda}y}\nonumber \\
 & =\sum_{i=0}^{w_{1}-1}\left(-1\right)^{i}\sum_{k=0}^{\infty}\mathcal{E}_{k}^{\left(m\right)}\left(\left.w_{2}x+\frac{w_{2}}{w_{1}}i\right|\frac{\lambda}{w_{1}}\right)w_{1}^{k}\frac{t^{k}}{k!}\nonumber \\
 & \relphantom =\times\sum_{l=0}^{\infty}\mathcal{E}_{l}^{\left(m-1\right)}\left(\left.w_{1}y\right|\frac{\lambda}{w_{2}}\right)w_{2}^{l}\frac{t^{l}}{l!}\nonumber \\
 & =\sum_{n=0}^{\infty}\sum_{k=0}^{n}\binom{n}{k}\sum_{i=0}^{w_{1}-1}\left(-1\right)^{i}\mathcal{E}_{k}^{\left(m\right)}\left(\left.w_{2}x+\frac{w_{2}}{w_{1}}i\right|\frac{\lambda}{w_{1}}\right)w_{1}^{k}\nonumber\\
 & \relphantom =\times\mathcal{E}_{n-k}^{\left(m-1\right)}\left(\left.w_{1}y\right|\frac{\lambda}{w_{2}}\right)w_{2}^{n-k}\frac{t^{n}}{n!}.\nonumber 
\end{align}

On the other hand, by the symmetric properties of $K^{\left(m\right)}\left(w_{1},w_{2}\mid\lambda\right)$
in $w_{1}$ and $w_{2}$, we get 
\begin{align}
 & K^{\left(m\right)}\left(w_{1},w_{2}\mid\lambda\right)\label{eq:26}\\
 & =\sum_{n=0}^{\infty}\sum_{k=0}^{n}\binom{n}{k}w_{2}^{k}w_{1}^{n-k}\mathcal{E}_{n-k}^{\left(m-1\right)}\left(\left.w_{2}y\right|\frac{\lambda}{w_{1}}\right)\nonumber\\
 & \relphantom =\times
\sum_{i=0}^{w_{2}-1}\left(-1\right)^{i}\mathcal{E}_{k}^{\left(m\right)}\left(\left.w_{1}x+\frac{w_{1}}{w_{2}}i\right|\frac{\lambda}{w_{2}}\right)\frac{t^{n}}{n!}.\nonumber 
\end{align}

Therefore, by comparing the coefficients on the both sides of (\ref{eq:25})
and (\ref{eq:26}), we obtain the following theorem.
\begin{thm}
\label{thm:4} For $w_{1},w_{2}\in\mathbb{N}$ with $w_{1}\equiv1,w_{2}\equiv1\pmod{2}$
and $n\ge0$ and $m\ge1$, we have 
\begin{align*}
 & \sum_{k=0}^{n}\binom{n}{k}w_{1}^{k}w_{2}^{n-k}\mathcal{E}_{n-k}^{\left(m-1\right)}\left(\left.w_{1}y\right|\frac{\lambda}{w_{2}}\right)\sum_{i=0}^{w_{1}-1}\left(-1\right)^{i}\mathcal{E}_{k}^{\left(m\right)}\left(\left.w_{2}x+\frac{w_{2}}{w_{1}}i\right|\frac{\lambda}{w_{1}}\right)\\
 & =\sum_{k=0}^{n}\binom{n}{k}w_{2}^{k}w_{1}^{n-k}\mathcal{E}_{n-k}^{\left(m-1\right)}\left(\left.w_{2}y\right|\frac{\lambda}{w_{1}}\right)\sum_{i=0}^{w_{2}-1}\left(-1\right)^{i}\mathcal{E}_{k}^{\left(m\right)}\left(\left.w_{1}x+\frac{w_{1}}{w_{2}}i\right|\frac{\lambda}{w_{2}}\right).
\end{align*}

\end{thm}
Let $y=0$ and $m=1$ in Theorem \ref{thm:4}. Then we have the following
corollary.
\begin{cor}
\label{cor:5} For $w_{1},w_{2}\in\mathbb{N}$ with $w_{1}\equiv1,w_{2}\equiv1\pmod{2}$,
$n\ge0$, we have 
\begin{align*}
 & w_{1}^{n}\sum_{i=0}^{w_{1}-1}\left(-1\right)^{i}\mathcal{E}_{n}\left(\left.w_{2}x+\frac{w_{2}}{w_{1}}i\right|\frac{\lambda}{w_{1}}\right)\\
 & =w_{2}^{n}\sum_{i=0}^{w_{2}-1}\left(-1\right)^{i}\mathcal{E}_{n}\left(\left.w_{1}x+\frac{w_{1}}{w_{2}}i\right|\frac{\lambda}{w_{2}}\right).
\end{align*}

\end{cor}
Let us take $w_{2}=1$ in Corollary \ref{cor:5}. Then we have 
\[
w_{1}^{n}\sum_{i=0}^{w_{1}-1}\left(-1\right)^{i}\mathcal{E}_{n}\left(\left.x+\frac{1}{w_{1}}i\right|\frac{\lambda}{w_{1}}\right)=\mathcal{E}_{n}\left(w_{1}x\mid\lambda\right).
\]


\bibliographystyle{spmpsci}      

\begin{thebibliography}{10}
\providecommand{\url}[1]{{#1}}
\providecommand{\urlprefix}{URL }
\expandafter\ifx\csname urlstyle\endcsname\relax
  \providecommand{\doi}[1]{DOI~\discretionary{}{}{}#1}\else
  \providecommand{\doi}{DOI~\discretionary{}{}{}\begingroup
  \urlstyle{rm}\Url}\fi

\bibitem{key-1}
Araci, S., Acikgoz, M.: A note on the {F}robenius-{E}uler numbers and
  polynomials associated with {B}ernstein polynomials.
\newblock Adv. Stud. Contemp. Math. (Kyungshang) \textbf{22}(3), 399--406
  (2012)

\bibitem{key-2}
Araci, S., Bagdasaryan, A., {\"O}zel, C., Srivastava, H.M.: New symmetric
  identities involving {$q$}-zeta type functions.
\newblock Appl. Math. Inf. Sci. \textbf{8}(6), 2803--2808 (2014).
\newblock \doi{10.12785/amis/080616}.
\newblock \urlprefix\url{http://dx.doi.org/10.12785/amis/080616}

\bibitem{key-3}
Bayad, A., Chikhi, J.: Apostol-{E}uler polynomials and asymptotics for negative
  binomial reciprocals.
\newblock Adv. Stud. Contemp. Math. (Kyungshang) \textbf{24}(1), 33--37 (2014)

\bibitem{key-4}
Bayad, A., Kim, T.: Identities involving values of {B}ernstein,
  {$q$}-{B}ernoulli, and {$q$}-{E}uler polynomials.
\newblock Russ. J. Math. Phys. \textbf{18}(2), 133--143 (2011).
\newblock \doi{10.1134/S1061920811020014}.
\newblock \urlprefix\url{http://dx.doi.org/10.1134/S1061920811020014}

\bibitem{key-5}
Carlitz, L.: Degenerate {S}tirling, {B}ernoulli and {E}ulerian numbers.
\newblock Utilitas Math. \textbf{15}, 51--88 (1979)

\bibitem{key-6}
Duran, U., Acikgoz, M., Araci, S.: Symmetric identities involving weighed
  $q$-{G}enocchi polynomials under $s_4$.
\newblock Proc. Jangjeon Math. \textbf{18}(4), (in press) (2015)

\bibitem{key-7}
He, Y.: Symmetric identities for {C}arlitz's {$q$}-{B}ernoulli numbers and
  polynomials.
\newblock Adv. Difference Equ. pp. 2013:246, 10 (2013).
\newblock \doi{10.1186/1687-1847-2013-246}.
\newblock \urlprefix\url{http://dx.doi.org/10.1186/1687-1847-2013-246}

\bibitem{key-8}
Kim, D.S., Kim, T.: Some identities of degenerate {E}uler polynomials arising
  from {$p$}-adic fermionic integrals on {$\Bbb{Z}_p$}.
\newblock Integral Transforms Spec. Funct. \textbf{26}(4), 295--302 (2015).
\newblock \doi{10.1080/10652469.2014.1002497}.
\newblock \urlprefix\url{http://dx.doi.org/10.1080/10652469.2014.1002497}

\bibitem{key-9}
Kim, D.S., Lee, N., Na, J., Park, K.H.: Identities of symmetry for higher-order
  {E}uler polynomials in three variables ({I}).
\newblock Adv. Stud. Contemp. Math. (Kyungshang) \textbf{22}(1), 51--74 (2012)

\bibitem{key-10}
Kim, D.S., Lee, N., Na, J., Park, K.H.: Abundant symmetry for higher-order
  {B}ernoulli polynomials ({I}).
\newblock Adv. Stud. Contemp. Math. (Kyungshang) \textbf{23}(3), 461--482
  (2013)

\bibitem{key-11}
Kim, T.: Symmetry of power sum polynomials and multivariate fermionic
  {$p$}-adic invariant integral on {$\Bbb Z_p$}.
\newblock Russ. J. Math. Phys. \textbf{16}(1), 93--96 (2009).
\newblock \doi{10.1134/S1061920809010063}.
\newblock \urlprefix\url{http://dx.doi.org/10.1134/S1061920809010063}

\bibitem{key-12}
Kim, T.: Symmetry properties of the generalized higher-order {E}uler
  polynomials.
\newblock Proc. Jangjeon Math. Soc. \textbf{13}(1), 13--16 (2010)

\bibitem{key-13}
Kim, T., Dolgy, D.V., Jang, Y.S., Seo, J.J.: A note on symmetric identities for
  the generalized {$q$}-{E}uler polynomials of the second kind.
\newblock Proc. Jangjeon Math. Soc. \textbf{17}(3), 375--381 (2014)

\bibitem{key-14}
Kim, T., Kim, D.S., Dolgy, D.V.: Degenerate {$q$}-{E}uler polynomials.
\newblock Adv. Difference Equ. p. 2015:246 (2015).
\newblock \doi{10.1186/s13662-015-0563-y}.
\newblock \urlprefix\url{http://dx.doi.org/10.1186/s13662-015-0563-y}

\bibitem{key-15}
Kim, Y.H., Hwang, K.W.: Symmetry of power sum and twisted {B}ernoulli
  polynomials.
\newblock Adv. Stud. Contemp. Math. (Kyungshang) \textbf{18}(2), 127--133
  (2009)

\bibitem{key-16}
Moon, E.J., Rim, S.H., Jin, J.H., Lee, S.J.: On the symmetric properties of
  higher-order twisted {$q$}-{E}uler numbers and polynomials.
\newblock Adv. Difference Equ. pp. Art. ID 765,259, 8 (2010)

\bibitem{key-17}
Rim, S.H., Jeong, J.H., Lee, S.J., Moon, E.J., Jin, J.H.: On the symmetric
  properties for the generalized twisted {G}enocchi polynomials.
\newblock Ars Combin. \textbf{105}, 267--272 (2012)

\bibitem{key-18}
Ryoo, C.S.: A note on the weighted {$q$}-{E}uler numbers and polynomials.
\newblock Adv. Stud. Contemp. Math. (Kyungshang) \textbf{21}(1), 47--54 (2011)

\bibitem{key-19}
{\c{S}}en, E.: Theorems on {A}postol-{E}uler polynomials of higher order
  arising from {E}uler basis.
\newblock Adv. Stud. Contemp. Math. (Kyungshang) \textbf{23}(2), 337--345
  (2013)

\end{thebibliography}


\end{document}